\documentclass [11pt]{article}
\usepackage{amsmath, amsthm, amssymb}
\usepackage{graphicx}
\newfont{\bbb} {msbm10}
\newcommand{\R}{\Bbb{R}}
\newcommand{\bS}{\Bbb{S}}
\newcommand{\bB}{\Bbb{B}}
\newcommand{\T}{\Bbb{T}}

\newcommand{\cT}{{\cal{T}}}

\newcommand{\sbs}{\subset}
\newcommand{\ra}{\rightarrow}

\newcommand{\bg}{{\bar{g}}}

\newcommand{\barB}{{\bar{B}}}

\newcommand{\hg}{{\hat{g}}}

\newcommand{\hH}{{\hat{H}}}

\newcommand{\p}{\partial}

\newcommand{\cH}{{\cal{H}}}

\newcommand{\cA}{{\cal{A}}}

\newcommand{\cW}{{\cal{W}}}

\newcommand{\B}{\Bbb{B}}

\newcommand{\HH}{\Bbb{H}}

\newcommand{\ssl}{{_\lambda}}

\newcommand{\0}[1]{_{_{#1}}}

\usepackage[all]{xy}

\begin{document}

\title{Deforming an $\epsilon$-Close to Hyperbolic Metric to a Hyperbolic Metric}
\author{Pedro Ontaneda\thanks{The author was
partially supported by a NSF grant.}}
\date{}

\maketitle

\begin{abstract} We show how to deform a metric of the form
$g=g_r+dr^2$ to a hyperbolic metric $\cH g$, 
for $r$ less than some fixed $\lambda$.
We study to what
extent the {\it hyperbolically forced metric} $\cH g$ is {\it close to being
hyperbolic}, if we assume $g$ to be close to hyperbolic.
We also deal with a one-parameter version of this problem.

The results in this paper are used in the problem of smoothing
Charney-Davis strict hyperbolizations \cite{ChD}, \cite{O}.

\end{abstract}

\noindent {\bf \large  Section 1. Introduction.}

First we introduce some notation.
Let $(M^n, g)$ be a complete Riemannian manifold {\it with center $o\in M$}, that is, 
the exponential map $exp_o:T_oM\ra M$ is a diffeomorphism.
In this case we
can write the metric $g$ on $M-\{ o\}=\bS^{n-1}\times\R^+$ as \,\,$g=g_r+dr^2$,
where $r$ is the distance to $o$. The open ball of radius $r$ in $M$, centered at $o$, will 
be denoted by $B_r=\B_r(M)$, and the closed ball by $\barB_r$.\vspace{.1in}

Let $M$ have center $o$ and metric $g=g\0{r}+dr^2$. In general, the metric $g$ is not a warped
product (i.e. $g_r=f(r)g'$ for some $g'$ not depending on $r$). In \cite{O2} we gave
a geometric process called {\it warp forcing}
whose input is a metric of the form
$g\0{r}+dr^2$ on $M$, and output is the {\it warped forced } metric $\cW\0{r\0{0}}g$ on $M$
(see Section 4 for more details). 
Here $r\0{0}>0$ is a fixed number. 
The metric $\cW\0{r\0{0}}g$ has the following property 

\noindent {\bf (1.1)}$\hspace{1.5in}\cW\0{r\0{0}} g\,=\,\left\{ \begin{array}{lllll}
\sinh^2r\,\,\hat{g}\0{r\0{0}}\, +\, dr^2&& {\mbox{on}}& & \barB_{r\0{0}}\\
g&&{\mbox{outside}}& & B_{r\0{0}+\frac{1}{2}}
\end{array}\right.$

\noindent where \vspace{.1in}

\noindent {\bf (1.2)}$\hspace{2.2in}\hat{g}\0{r\0{0}}=
\big(  \frac{1}{\sinh^2r\0{0}}\big) \,g\0{r\0{0}}$\vspace{.1in}

Hence warp forcing does not change
the metric outside $B_{r\0{0}+\frac{1}{2}}$ while making it a 
warped product inside $\barB_{r\0{0}}$. \vspace{.1in}

\noindent {\bf Remark 1.3.}  We can think of the metric $g\0{r\0{0}}$ in
as being obtained from $g=g_r+dr^2$ by ``cutting" $g$ along
the sphere of radius $r\0{0}$, so we call
$g\0{r\0{0}}$ the {\it spherical cut of
$g$ at $r\0{0}$.}
Also we call the metric $\hat{g}\0{r\0{0}}$ in (1.2) the {\it normalized spherical
cut of $g$ at $r\0{0}$} (see Section 4).
\vspace{.1in}

On the other hand in \cite{O1} we introduced another geometric construction,
the {\it two variable warping trick} (see Section 3
for more details). Its input is a warped product $g=\sinh^2r\,\, h +dr^2$
on $\bS^{n-1}\times\R^+$, and its output is the metric $\cT\0{a,d}g$, where $a,d>0$ are two
parameters. The metric $\cT\0{a,d}g$ has the following property. 

\noindent {\bf (1.4)}$\hspace{1.5in}\cT\0{a,d}g\,=\,\left\{ \begin{array}{lllll}
\sinh^2r\,\,\sigma\0{\bS^{n-1}}\, +\, dr^2&& {\mbox{on}}& & \barB_{a}\\
g&&{\mbox{outside}}& & B_{a+\frac{d}{2}}
\end{array}\right.$

\noindent Here $\sigma\0{\bS^{n-1}}$ is the canonical round metric on the sphere.
Hence, the two variable warping trick
changes a warped product  $g$ inside the ball $\barB_{a+\frac{d}{2}}$ 
making it (radially) hyperbolic on the smaller ball $\barB_a$.
The warped product $g$ does not change outside 
$B_{a+\frac{d}{2}}$.\vspace{.1in}

\noindent {\bf Remark 1.5} To be able to define $\cT\0{a,d}g$ 
the metric $g$ does not need to be a warped product
everywhere. It only needs to be a warped product in the ball $\barB_{a+\frac{d}{2}}$.
Also to define $\cW\0{r\0{0}} g$, the metric $g$ needs to be defined only
outside $B_{r_0}$. See sections 3 and 4.\vspace{.1in}

In this paper we consider the composition of these two geometric constructions.
Given $g=g\0{r}+dr^2$,  $d>0$, $r\0{0}$ with $r\0{0}>d$,  we define
the  metric $\cH\0{r\0{0},d}\,\big( g)$ in the following way.
First warp-force the metric $g$, i.e take $\cW\0{r\0{0}} g$. Since $\cW\0{r\0{0}} g$
is a warped product on $B\0{r\0{0}}$ (see (1.1)) we can use the two variable warping trick
(see Remark 1.5) and define\vspace{.1in}

\noindent {\bf (1.6)}$\hspace{2.2in}\cH\0{r\0{0},d}\,\, g\,=\, \cT_{_{(r\0{0}-d),d}}
\big( \,\, \cW\0{r\0{0}} g\,\,   \big)$\vspace{.1in}

\noindent  We call {\it hyperbolic forcing
(with cut at $r\0{0}$}) the process given by
the correspondance $g\mapsto \cH\0{r\0{0},d}\, g$.
It follows from (1.1) and (1.4) that

\noindent {\bf (1.7)}$\hspace{1.5in}\cH\0{r\0{0},d}g\,=\,\left\{ \begin{array}{lllll}
\sinh^2r\,\,\sigma\0{\bS^{n-1}}\, +\, dr^2&& {\mbox{on}}& & B\0{r\0{0}-d}\\
{\mbox{warped product}}&&{\mbox{on}}& &B\0{r\0{0}}- B\0{r\0{0}
-\frac{d}{2}}\\
g {\mbox{\,\,(no change)}} && {\mbox{outside}}&& B\0{r\0{0}
+\frac{1}{2}}
\end{array}\right.$\vspace{.1in}

\noindent Hence the metric $\cH\0{r\0{0},d}\, g$ is hyperbolic inside $\barB_{r\0{0}-d}$.
In Section 4 we give some properties of this metric.
\vspace{.1in}

In this paper we prove that if $g$ is in some sense
``close to being hyperbolic" then
the hyperbolically forced metric  $\cH\0{r\0{0},d}\,$ is also, to some extent,  close to 
hyperbolic. In the next paragraph
we explain what we mean by a metric being close to hyperbolic. 
Essentially it is a chart-by-chart concept.
\vspace{.1in}

Let $\bB\sbs\R^{n-1}$ be the open unit $(n-1)$-ball, 
with the flat metric $\sigma\0{\R^{n-1}}$.  Write  $I_\xi=(-1-\xi,1+\xi)$, $\xi>0$.
Our basic models are $\T_\xi=\bB\times I_\xi$,  with hyperbolic metric $\sigma=e^{2t}\sigma\0{\R^{n-1}}+dt^2$. 
The number $\xi$ is called the {\it excess} of $\T_\xi$.
Let $(M,g)$ be a Riemannian manifold and $S\sbs M$. We say that $S$ is $\epsilon$-{\it close to hyperbolic} if there is $\xi>0$ such that for every $p\in S$ there is 
an {\it $\epsilon$-close to hyperbolic 
chart with center $p$}, that is, there is a chart
$\phi :\T_\xi\ra M$, $\phi(0,0)=p$,  such that 
$|\phi^*g-\sigma|\0{C^2}<\epsilon$, where $|.|\0{C^2}$ is the $C^2$ norm (see
Section 2).  
The number $\xi$ is called the {\it excess } of the charts
(which  is fixed). \vspace{.1in}

\noindent {\bf Remark 1.8.}
Here is a simple but useful observation that will be used later.
Suppose two metrics $g_1$, $g_2$ on $M$ coincide on the ball
$B_{1+\xi}(p)$ of radius $1+\xi$ centered at $p$. Also suppose that $\phi$ is an $\epsilon$-close to hyperbolic chart centered at $p$ with excess $\xi$ for 
$(M,g_1)$. Then $\phi$ is also an $\epsilon$-close to hyperbolic chart centered at $p$ with excess $\xi$ for 
$(M,g_2)$. \vspace{.1in}

Let $(M,g)$ have center $o$ and $S\sbs M$. We say that $g$ is {\it radially $\epsilon$-close to hyperbolic on $S$ (with respect to $o$)} if, 
for every $p\in S$ there is an  $\epsilon$-close to hyperbolic 
chart $\phi$ with center $p$ and, in addition, the chart 
$\phi$ respects the product structure of $\T_\xi$
and $M-o=\bS^{n-1}\times\R^+$, that is $\phi(. , t)=(\phi\0{1}(.), t+a)$, 
where the constant 
 $a$ depends on 
$\phi$, and $\phi\0{1}$ is some function independent of $t$ (equivalently, $\phi_1$ is a chart on $M$). Here the ``radial" directions
are $(-1-\xi,1+\xi)$ and $\R^+$ in $\T_\xi$ and $M-o$, respectively.
\vspace{.1in}

\noindent{\bf Remark 1.9.} 
We prove in 4.14 of \cite{O1} that hyperbolic $n$-space
$\HH^n$ is radially $\epsilon$-close to hyperbolic (with respect to a fixed
center $o$ and with charts of excess $\xi$) outside the ball $\barB_{\sf a}(o)$,
where ${\sf a}$ is a constant that depends on $\epsilon$, $n$ and $\xi$.
But $\HH^n$ is not radially $\epsilon$-close to hyperbolic on
the ball $\barB_{\sf a}(o)$. In general we will use the notion
of radially $\epsilon$-close to hyperbolic metric
for large $t$.\vspace{.1in}

As mentioned above our first result says that if the metric $g=g_r+dr^2$
is radially $\epsilon$-close to hyperbolic then
the hyperbolically forced metric $\cH\0{r\0{0},d}\,g$
is radially $\eta$-close to hyperbolic, where $\eta$
depends on $\epsilon$, $r\0{0}$ and $d$. Recall that
$\cH\0{r\0{0},d}\,g$ is hyperbolic on $\barB_{r\0{0}-d}$
(see (1.7)).
This motivates us to deal with a natural and useful class of metrics.
These are metrics on $\R^n$ (or on a manifold with center)
that are already hyperbolic on the ball
$\barB_{a}=\barB_a(0)$ of radius $a$ centered at $0$, and are radially $\epsilon$-close 
to hyperbolic outside
$\barB_{a'}$
(here $a'$ is slightly less than $a$). Here is the detailed
definition. Let $M^n$ have center $o$ and
let $B_{a}=B_a(o)$ be the ball on $M$ of radius $a$ centered at $o$. We say
that a metric $h$ on $M$ is $(B_a,\epsilon)$-{\it close to hyperbolic with 
charts of excess $\xi$}, \, if 
\begin{enumerate}
\item[  (1)]  On $\barB_{a}-\{o\}=\bS^{n-1}\times (0,a)$
we have $h=\sinh^2r\,\,\sigma\0{\bS^{n-1}}+dr^2$. Hence $h$
is hyperbolic on $\barB_a$.
\item[(2)]  the metric $h$ is 
radially $\epsilon$-close to hyperbolic outside $\barB_{a-1-\xi}$
with charts of excess $\xi$.
\end{enumerate}

We will always assume $a>{\sf a}+1$,
where {\sf a} is as in 1.9 (this implies that $a-1-\xi>0$, see p. 579 of \cite{O1}). 
We have dropped the word ``radially" in the definition above
to simplify the notation; but it does appear in condition (2),
where ``radially" refers to the center of $B_a$.\vspace{.1in}

Metrics that are $(B_a,\epsilon)$-close to
hyperbolic are very useful, and are key objects in
\cite{O}. See also \cite{O1}, \cite{O3}.
Our next result answers the following question:
\vspace{.1in}

\noindent {\it To what extend is the hyperbolically forced metric 
$\cH\0{r\0{0},d}\, g$ close to hyperbolic, when
$g$ is close to hyperbolic?}\vspace{.1in}

Before we state our first result we need one more definition.
A family $g\0{t}$ of metrics on the sphere is {\it $c$-bounded} if each metric $g\0{t}$ 
(and its derivatives up to
order 2) are bounded by $c$ and $|det\,g\0{t}|_{C^0}>1/c$ (see Section 2 for more details). 
We are ready to state our 
first result. Recall that $\hat{g}\0{r\0{0}}$ is the normalized spherical
cut of $g$ at $r\0{0}$ (see 1.2). 
\vspace{.1in}

\noindent {\bf Theorem 1.} {\it Let $M^{n}$ have center $o$ and 
metric $g=g\0{r}+dr^2$ and\, $r\0{0},\, d >0$. Assume 
the normalized spherical cut $\hg\0{r\0{0}}$ is $c$-bounded.
If the metric $g$ is radially $\epsilon$-close to hyperbolic 
outside $\barB\0{r\0{0}-1-\xi}$
with charts of excess $\xi>1$, then the metric $\cH\0{r\0{0},d}\, g$ is 
$(B_{r\0{0}-d},\eta$)-close to hyperbolic
with charts of excess $\xi-1$ provided} \begin{center}$\eta\geq
C_1\,\Big(\,\frac{1}{d}\,+\,e^{-(r\0{0}-d)}\Big)\,+\, C_2
\,\epsilon$\end{center}
\noindent {\it  Here 
$C_1$ is a constant depending only on $n$, $\xi$,  $c$, and $C_2$ depends 
only on $\xi$. Also we are assuming $r\0{0}>d\geq 4+4\xi$.}
\vspace{.1in}

\newpage

\noindent {\bf  Remarks 1.10.}

\noindent {\bf 1.} An important point here is that by taking 
$r\0{0}$ and $d$ large the metric $\cH\0{r\0{0},d}\, g$ can be made 
$2C_2\epsilon$-close 
to hyperbolic.
How large we have to take $d$ and $r\0{0}$ depends on $c$, which is a $C^2$ bound
for the $\hat{g}\0{r\0{0}}$, the normalized spherical
cut of $g$ at $r\0{0}$.

\noindent {\bf 2.} 
We give a formula for $C_1=C_1(c,n,\xi)$ at the end of Section 4.  
We can take $C_2=e^{27+12\xi}$.
\vspace{.1in}

Next we give a one-parameter version of hyperbolic forcing
with cut at $r\0{0}$, with the
variable $r\0{0}\ra \infty$ (we change notation and use $\lambda$ 
instead of $r\0{0}$
to express the fact that this number now varies). This family 
version of hyperbolic forcing
is an important ingredient in \cite{O}, where it is used
to smooth the singularities of a
Charney-Davis strict hyperbolization \cite{ChD}.
\vspace{.1in}

Let $M^n$ have center $o$. As before write $M-\{o\}=\bS^{n-1}\times\R^+$. 
Recall that to 
apply hyperbolic forcing, with cut at $r\0{0}$, to a metric $g=g\0{r}+dr^2$, 
the metric $g$
needs to be defined only for $r>r'$, where $r'<r\0{0}$ (see Remark 1.5). 
\vspace{.1in}

Fix $\xi$. We now consider families of metrics 
$\{ g_\ssl\}_{\lambda> \lambda_0}$ of the form $g_\ssl=\big(g_\ssl\big)_r+dr^2$. 
Here $\lambda_0>1+\xi$.
  We will assume that
$\big(g_\ssl\big)_r$ defined for (at least) $r>\lambda-1-\xi$. 
We call such a family
a $\odot${\it-family of metrics}. We now apply hyperbolic
forcing to each $g_\ssl$ with cut a $r=\lambda$,
that is we consider  $\cH\0{\lambda,d}g\0{\lambda}$.
Note that to obtain $\cH\0{\lambda,d}g\0{\lambda}$ from $g_\ssl$ we are cutting at 
$\lambda$, that is we are cutting $g_\ssl$ ``along"
the sphere of radius $\lambda$.\vspace{.1in}

We want to give a one-paramenter version of
Theorem 1, that is, a version for a $\odot$-family $\{ g_\ssl\}$. 
Since the constant
$C_1$ in Theorem 1 depends on the bound $c$ there is no uniform $C_1$
that would work for {\sf every}
$g_\ssl$. This problem motivates the following
definition.\vspace{.1in}

We say that the  $\odot$-family $\{g_\ssl\}$ has {\it cut limit at $b$} 
if there is a 
$C^2$ metric $\hg_{_{\infty}}^b$ on $\bS^{n-1}$ such that 
\begin{center}$
{\widehat{\big(g_\ssl\big)}}_{_{\lambda+b}} \,\,\,\,\stackrel{C^2}
{\longrightarrow}\,\,\,\, \hg_{_{\infty}}^b
$\end{center}
\noindent 
as $\lambda\ra\infty$. Recall that the ``hat" notation means ``normalized spherical
cut" (see 1.2). One  more definition. We say that
an $\odot$-family $\{g_\ssl\}$  is radially $\epsilon$-close to hyperbolic
with charts of excess $\xi$, if each
$g\0{\lambda}$ is radially $\epsilon$-close to hyperbolic outside 
$\barB_{\lambda-1-\xi}$
with charts of excess $\xi$. The following is a corollary of Theorem 1, 
and is proved in
Section 6.\vspace{.1in}

\noindent {\bf Theorem 2.} {\it Let $M$ have center $o$, and $\{g\0{\lambda}\}$
an $\odot$-family on $M$. Assume that $\{g\0{\lambda}\}$ has cut limits at $b=0$.
Then there is $c$ such
that the following holds. 
Suppose $\{ g\0{\lambda}\}$ is radially $\epsilon$-close to hyperbolic with charts 
of excess $\xi>1$. Then, for every $\lambda$,   $\cH\0{\lambda,d}g\0{\lambda}$ is  
$(B_{\lambda-d},\eta)$-close to hyperbolic
with charts of excess $\xi-1$, where } \begin{center}$\eta\geq
C_1\,\Big(\,e^{-(\lambda-d)}+\,\frac{1}{d}\,\Big)\,+\,
C_2\,\epsilon$\end{center}
\noindent {\it
Here $C_2=C_2(\xi)$ and $C_1=C_1(c, n,\xi)$ are as in Theorem 1.}\\

\noindent {\bf Corollary.} {\it Let $M$ have center $o$, 
$\{g\0{\lambda}\}$
an $\odot$-family on $M$, and $\epsilon'>0$. Assume that $\{g\0{\lambda}\}$ 
has cut limits at $b=0$.
If $\{ g\0{\lambda}\}$ is radially $\epsilon$-close to hyperbolic with charts 
of excess $\xi>1$, then, for every $\lambda$,  $\cH\0{\lambda,d}g\0{\lambda}$ is 
$(B_{\lambda-d},\epsilon'+C_2\epsilon)$-close to hyperbolic,
with charts of excess $\xi-1$, provided
\begin{enumerate}
\item[(i)] $\lambda-d>\ln(\frac{2C_1}{\epsilon'})$  
\item[(ii)] $d\geq \frac{2C_1}{\epsilon'}$.
\end{enumerate}

\noindent Here $C_1$ and $C_2$ are as in Theorem 2.}\vspace{.1in}

Note that we can take $\epsilon'$ as small as we want hence $\epsilon'+C_2\epsilon$
as close as $C_2\epsilon$ as we desire, provided we choose $d$ and $\lambda$
sufficiently large. How large depending on $\epsilon'$
and $c$.\vspace{.1in}

To deduce Corollary from Theorem 2 just note that 
(i) and (ii) in Corollary imply $C_1e^{-(\lambda-d)}<\epsilon'/2$
and $C_1\frac{1}{d}<\epsilon'/2$.\vspace{.1in}

The results in this paper are used to construct negatively curved Riemannian smoothings of Charney-Davis strict
hyperbolizations of manifolds \cite{ChD}, \cite{O}. In the next paragraph we give an idea how the objects and results in this paper is used in \cite{O}.
\vspace{.1in}

In the same way that a cubical complex is made of basic pieces (the cubes
$\square^k$),
the hyperbolization $h(K)$ of a cubical complex $K$ is 
also made of
basic pieces: pre-fixed hyperbolization pieces $X^k$. Indeed one begins with a cubical complex $K$ and replaces each cube of dimension $k$
by the hyperbolization piece of the same dimension. Cube complexes have
a piecewise flat metric induced from the flat geometry of the cubes.
Likewise the Charney-Davis hyperbolizations have a piecewise hyperbolic
structure because the Charney-Davis hyperbolization pieces
are hyperbolic manifolds (compact, with boundary and corners).
To see how singularities appear one can first think about the manifold 2-dimensional 
cube case. If $K^2$ is a 2-dimensional manifold cube complex then
its piecewise flat metric is Riemannian outside the vertices. A vertex is
a singularity if and only if the vertex does not meet exactly four cubes.
The picture is exactly the same for $h(K^2)$. 
These point singularities in $h(K^2)$ can be smoothed out
easily using warping methods.
In higher dimensions
the singularities of $K^n$ and $h(K)$ appear in (possibly the whole of)
the codimension 2 skeletons $K^{(n-2)}$ and $h(K^{(n-2)})$, respectively.
In \cite{O} the idea of smoothing the piecewise hyperbolic metric on
$h(K)$ is to do it inductively down the dimension of the skeleta.
One begins with the $(n-2)$-dimensional pieces $X^{n-2}$. Transversally
to each $X^{n-2}$ (that is, on the union of geodesic segments emanating
perpendicularly to $X^{n-2}$, from a fixed point in $X^{n-2}$) one has
essentially the 2-dimensional picture mentioned above. Once we
solve this transversal problem we extend this transversal smoothing by taking a warp product
with $X^{n-2}$; we called this product method {\it hyperbolic extension}
\cite{O3}.
This gives a smoothing on a (tubular) neighborhood of the piece $X^{n-2}$.
Caveat: we do not want to actually have a smoothing on a neighborhood
of the whole of $X^{n-2}$, since we will certainly have 
matching problems for different $X^{n-2}$ meeting on a common
$X^{n-3}$; so we only want a smoothing on a neighborhood of the $Z^{n-2}$,
where $Z^{n-2}\sbs X^{n-2}$ is just a bit ``smaller" than $X^{n-2}$,
so that the neighborhoods of the $Z^{n-2}$ are all disjoint. Next 
step is to smooth around the $X^{n-3}$ (or, specifically the $Z^{n-3}$).
The metric is already smooth outside a neighborhood of the $(n-3)$-skeleton. Transversally to each $X^{n-3}$ we
have a 3 dimensional problem.
(It helps to have a 3 dimensional picture in mind, like in dimension 2).
It happens that if we did things with care in the first step (around the
$Z^{n-2}$) the metric in the 3 dimensional transversal problem
is radially $\epsilon$-close to hyperbolic outside some large ball B.
At this point we use the method of {\it hyperbolic forcing} introduced in
this paper to
 to extend the metric to a Riemannian metric
on the ball B, getting rid, in this way, of the transverse
singularity. The resulting
metric is still radially $\eta$-close to hyperbolic, with an $\eta$
that can be controlled. Once the transversal 3 dimensional problem
is solved we extend this smoothing to neighborhoods of the $Z^{n-3}$
using hyperbolic extension. Next we do the same for the $Z^{n-4}$
and so on. About the excess: since hyperbolic forcing reduces the excess
by 1, one begins with a large excess at codimension 2, so that when
we arrive at codimension $n$ one still has positive excess;
therefore in the Theorem above one
should think of the $\xi$ as fixed, while the $r\0{0}$ as being
as large as wanted, $\epsilon$ as small as desired, and the set $S$ as the complement of the
ball of radius $r\0{0}-1-\xi$.\vspace{.1in}

In Section 2 we give a bit more details about the $C^k$-norm
used here. In Section 3 we deal briefly with the two variable warping trick. In 
Section
4 we do the same with the warp forcing process. In Section 5 we study
hyperbolic forcing and prove Theorem 1. In Section 6 we deal with 
$\odot$-families and prove Theorem 2.
\vspace{.3in}

\noindent {\bf \large  Section 2. The $C^k$-Norm.}\\ Let $A\sbs \R^n$ be an open set.
Let $|.|\0{C^2(A)}$ denote the uniform $C^2$-norm of $\R^l$-valued functions on $A$, i.e. if $f=(f\0{1},...,f\0{l}):A\ra \R^l$,
then $|f|\0{C^2(A)}=sup\0{z\in A,
\,\,1\leq i\leq l,\,\,1\leq j,k\leq n}\{ |f\0{i}(z)|, |\p\0{j}f\0{i}(z)|, |\p\0{j,k}f\0{i}(z)|\}$. Sometimes we will write $|.|\0{C^2}=|.|\0{C^2(A)}$ when the context is clear.
Given a Riemannian metric $g$ on $A$, the number $|g|\0{C^2(A)}$ is computed considering $g$ as the $\R^{n^2}$-valued function $z\mapsto (g_{ij}(z))$ where, as usual,
$g\0{ij}=g(e_i,e_j)$, and the $e_i$'s are the canonical vectors in $\R^{n}$. 
\vspace{.1in}

The $C^2$-norm $|.|\0{C^2}$ mentioned in the  
definition of an  $\epsilon$-close to hyperbolic Riemannian manifold in the Introduction is $|.|\0{C^2}=|.|\0{C^2(\T_\xi)}$. 
If $(M,g)$ is $\epsilon$-close to hyperbolic (or radially $\epsilon$-close 
to hyperbolic) we will also say that the metric $g$ is $\epsilon$-close to hyperbolic (or radially $\epsilon$-close to hyperbolic).
\vspace{.1in}

Let $c>1$.
A metric $g$ on a compact manifold $M$ is {\it c-bounded} if $|g|\0{C^2}< c$ and 
$|\, det \,g\,|_{C^0} > 1/c$. 
A set of metrics $\{ g\0{\lambda}\}$ on the compact manifold $M$ is 
{\it c-bounded} if every $g\0{\lambda}$
is $c$-bounded.  \vspace{.1in}

\noindent {\bf Remarks 2.1.}

\noindent {\bf 1.} Here the uniform $C^k$-norm $| .|\0{C^k}$ is taken with respect 
to a fixed locally finite atlas $\cA$ with ``extendable" charts, i.e.
charts that can be extended to the (compact) closure of their domains.

\noindent {\bf 2.} If $\{ g\0{t}\}_{t\in I}$ is a smooth family and $I\sbs \R$ 
is compact then, by compactness of $M$ the family
 $\{ g\0{t}\}$ is $c$-bounded for some $c$. 
\vspace{.3in}

\noindent {\bf \large 3.  The Two Variable Warping Trick.} 

Here we review the two variable warping trick. For more details see \cite{O1}.
Let $h$ be a metric on the $(n-1)$-sphere $\bS^{n-1}$ and consider the warped product
$g=\sinh^2t\, h +dt^2$ on $\bS^{n-1}\times \R^+$. We fix a function $\rho:\R\ra[0,1]$ 
with $\rho (t)=0$ for $t\leq 0$, $\rho(t)=1$ for $t\geq 1$, and $\rho$ constant near 0 and 1.  
Given positive numbers $a$ and $d$ define $\rho\0{a,d}(t)=\rho(2\,\frac{t-a}{d})$. 
Also fix an atlas $\cA\0{\bS^n}$ on $\bS^{n-1}$ as before (see Remark 2.1(1)).
All norms and boundedness constants will be taken with respect to
this atlas.\vspace{.1in}

As before let $\sigma\0{\bS^{n-1}}$ be the round metric on $\bS^{n-1}$.
Write \vspace{.1in}

\noindent {\bf  (3.1)}\hspace{2in}$h\0{t}=\big(\,1-\rho\0{a,d}(t)\, \big) 
\sigma\0{\bS^{n-1}} +\rho\0{a,d}(t) \,h$\vspace{.1in}

\noindent  and define the metric\vspace{.1in}

\noindent {\bf  (3.2)}\hspace{2in}$
\cT_{_{a,d}}\, g\, =\, \sinh^2\, t\,\,h\0{t}+dt^2
$\vspace{.1in}

The process $g\mapsto \cT\0{a,d}g$ is the {\it two variable warping trick}. We have that 
$\cT\0{a,d}g$ satisfies property (1.4) in the Introduction.\vspace{.1in}

\noindent {\bf Remark 3.3.} Note that to define $\cT\0{a,d}g$ we need $g$ to be
a warped product only for $r$ near $a+d$. 
\vspace{.1in}

The following  is Corollary 5.6 in \cite{O1} \vspace{.1in}

\noindent {\bf Theorem 3.4.} {\it Let the metric $h$ on $\bS^{n-1}$ be $c$-bounded. 
Write $g=\sinh^2t\,\,h
+dt^2$.  Then the metric $\cT_{_{a,d}}\, g\,$ is  $(B_a,\epsilon)$-close to 
hyperbolic with charts of excess $\xi$, provided} 
\begin{center}$\,C\,\Big(e^{-a}+\frac{1}{d}\,\Big)\leq\epsilon$
\end{center}
\noindent {\it Here $C$ is a constant depending on
$c$, $n$ and $\xi$.} \vspace{.1in}

\noindent For an explicit formula of $C$ see
\cite{O1} (the constant $C$ here is denoted by
$C_2$ in \cite{O1}). The  Corollary in the Introduction of
\cite{O1} is obtained 
from 4.4 of \cite{O1} (Theorem 3.4 here)
by taking $\xi=0$. 

\vspace{.1in}


\vspace{.3in}

\noindent {\bf \large  Section 4. Spherical Cuts and Warp Forcing.}

As in the Introduction, let $(M^n, g)$ be a complete Riemannian manifold with 
center $o\in M$. Recall that we
can write the metric on $M-\{ o\}=\bS^{n-1}\times\R^+$ as \,\,$g=g_r+dr^2$. 
\vspace{.1in}


Fix $r\0{0}>0$.
We define the warped by $\sinh$ metric $\bg\0{r\0{0}}$ by:\vspace{.1in}

\noindent {\bf (4.1)}\hspace{2in}$\bg \0{r\0{0}}\,=\,\sinh^2r\,\,\,\hg\0{r\0{0}}\,+\, dr^2     
$\vspace{.1in}

\noindent where $\hat{g}\0r\, =\, \big(\frac{1}{\sinh^2r}\big)\, 
g_r$ is the normalized spherical cut of $g$ at $r$ (see 1.3).
We now force the metric $g$ to be equal to $\bg\0{r\0{0}}$ on  
the ball $\barB_{r\0{0}}=\barB_{r\0{0}}(M)$ and stay equal to $g$ outside $B_{r\0{0} +\frac{1}{2}}$.
For this we define the {\it warped forced } (on $B_{r\0{0}}$) metric as:\vspace{.1in}

\noindent {\bf (4.2)}\hspace{1.9in}
$\cW\0{r\0{0}}\, g\, =\, (1-\rho\0{r\0{0}}) \, \bg\0{r\0{0}}\, +\, \rho\0{r\0{0}} \, g
$\vspace{.1in}

\noindent where $\rho\0{r\0{0}}(t)=\rho(2t-2r\0{0})$, with $\rho$ as in Section 3,
that is,  $\rho:\R\ra[0,1]$ 
with $\rho (t)=0$ for $t\leq 0$, $\rho(t)=1$ for $t\geq 1$, and $\rho$ constant near 0 and 1.  
Hence we have that $\cW\0{r\0{0}} g$
satisfies (1.1) in the Introduction.\vspace{.1in}

Therefore the warp forced metric $\cW\0{r\0{0}} g$ is a warped product on $\barB_{r\0{0}}$.
We call the process $g\mapsto\cW g$ {\it warp forcing}.
The next result is the Main Theorem of \cite{O2}. It states that if $g$ is 
{\it $\epsilon$-close to a hyperbolic metric} then
the warp forced metric  $\cW\0{r\0{0}} g$ is also close to hyperbolic. \vspace{.1in}

\noindent {\bf Theorem 4.3.} {\it  Let $(M,g)$ have center $o$, and $S\sbs M$. 
If $g$ is radially $\epsilon$-close to hyperbolic on $S$ 
with charts of excess $\xi$,
then $\cW\0{r\0{0}} g$ is radially $\eta$-hyperbolic on $S-\barB\0{r\0{0}-1-\xi}$ 
with charts of excess $\xi-1$, provided 
$\eta\geq e^{27+12\xi}\big(e^{-2r\0{0}}+\epsilon\big)$.}
\vspace{.1in}

\noindent{\bf Remarks 4.4.} 


\noindent {1.} Note that to define $\cW\0{r\0{0}} g$ we only need
$g\0{r}$  defined for $r\geq r\0{0}$. But to apply Theorem 4.3 we need 
$g\0{r}$ defined for $r\geq r\0{0}-1-\xi$.

\noindent {2.} Notice that $\cW\0{r\0{0}} g$ is a warped product no just on
$B\0{r\0{0}}$ but on $B\0{r\0{0}+\delta}$, for small $\delta>0$.\vspace{.3in}

\noindent {\bf \large  Section 5. Hyperbolic Forcing and Proof of Theorem 1.}

Let $(M^{n},g)$ have center $o$. As before we write $g=g\0{r}+dr^2$.
Let  $r\0{0}>d>0$. Recall that in the Introduction (see (1.6) 
we defined
the  metric $\cH\0{r\0{0},d}g$ by
$\cH\0{r\0{0},d}\,g\,=\, \cT_{_{(r\0{0}-d),d}}\big( \, \cW\0{r\0{0}} g\,\big)$, and called 
the process $g \mapsto \cH\0{r\0{0},d}g$  {\it hyperbolic forcing}. Note that $\cH\0{r\0{0},d}g$ satisfies
(1.7).\vspace{.1in}

Write $H=\cH\0{r\0{0},d}\,g\,$. 
Note that $H$ also has the form  $H=H_r+dr^2$.
We can explicitly describe $ H_r$:\vspace{.1in}

\noindent {\bf Proposition 5.1.} {\it We have}
\vspace{.1in}

\noindent \,\,\,\,\,\,\,\,\,\,\,\,\,{\small $
H_r\,\, =\,\,\left\{
\begin{array}{lll}
g_r&&   r\0{0} +\frac{1}{2}\leq r\\
\big(1-\rho\0{r\0{0}} (r)\big) \,\sinh^2r\,\,\, \hg\0{r\0{0}} \, +\,\rho\0{r\0{0}} (r) \,g\0{r}& 
&r\0{0} \leq r\leq r\0{0} +\frac{1}{2}\\
\sinh^2 r\, \bigg( \,\big( 1-  \rho_{_{(r\0{0} -d), d}}(r)\, \big)\, 
\sigma\0{\bS^{n-1}}\,+\, \rho_{_{(r\0{0} -d), d}}(r) \, \hg\0{r\0{0}}\,  \bigg)&&  r\0{0} -d 
\leq r\leq r\0{0}  \\
\sinh^2r\,\,\, \sigma\0{\bS^{n-1}}&&r\leq r\0{0}-d
\end{array}
\right.
$}\vspace{.1in}

\noindent {\it where the gluing functions $\rho\0{r\0{0}}$ and $ \rho_{_{(r\0{0} -d), d}} $ 
are defined in sections 4 and 3, respectively.}\vspace{.1in}

\noindent {\bf Proof.} We have four cases.

\noindent {\bf (i)}  For $r\geq r\0{0}+1/2$ we 
have that $H_r=g_r$. This follows from (1.7).

\noindent {\bf (ii)} For $r\0{0}\leq r\leq r\0{0}+1/2$
we have, by (1.4), that $H_r=\big(\cW_{r\0{0}}g\big)_r$. And by (4.2), (4.1), 
we get that $\big(\cW_{r\0{0}}g\big)_r=
\big(1-\rho\0{r\0{0}} (r)\big) \,\sinh^2r\,\,\, \hg\0{r\0{0}} \, 
+\,\rho\0{r\0{0}} (r) \,g\0{r}$.

\noindent {\bf (iii)} For $r\0{0}-d\leq r\leq r\0{0}$
we have, by (4.1) and (1.1), that $\big(\cW_{r\0{0}}g\big)_r=\sinh^2r\,\,\hg\0{r\0{0}}$. 
This together with (3.1) and (3.2) 
gives $H_r=\sinh^2r\,\, \big( \,\big( 1-  \rho_{_{(r\0{0} -d), d}}(r)\, \big)\, 
\sigma\0{\bS^{n-1}}\,+\, \rho_{_{(r\0{0} -d), d}}(r) \, \hg\0{r\0{0}}\,  \big)$.

\noindent {\bf (iv)} For $ r\leq r\0{0}-d$
we have, by (1.4), that $H_r=\sinh^2r\,\, \sigma\0{\bS^{n-1}}$. 
This proves the proposition.
\vspace{.2in}

\noindent {\bf Proposition 5.2.}  {\it The metric $H=\cH\0{r\0{0},d}\,g\,$ 
has the following properties.}

{\it 
\noindent (i) The metric $H$ is canonically hyperbolic on $\barB_{_{r\0{0}-d}}$, i.e. 
$H=\sinh^2r\,\,\sigma\0{\bS^n}+dr^2$ on
$\barB\0{r\0{0}-d}$.

\noindent (ii) We have that $g=H$ outside $B_{_{r\0{0}+\frac{1}{2}}}.$

\noindent (iii) The metric $H$ coincides with $\cW\0{r\0{0}} \big(g\0{r\0{0}}\big)$ 
outside $B\0{r\0{0}-\frac{d}{2}}$.

\noindent (iv) The metric $H$ coincides with 
$\cT_{_{(r\0{0}-d),d}} \bar{g}\0{r\0{0}}$
 on $\barB\0{r\0{0}}$.

\noindent (v) All the $g$-geodesic rays $r\mapsto ru$, $u\in \bS^{n-1}$, emanating 
from the center are geodesics of
$(M , H)$. Hence, the space  $(M,H)$ has center $o$. Moreover
the function $r$ (distance to the center $o$) is the same on the spaces $(M,g)$ 
and $(M,H)$.}\vspace{.1in}

\noindent {\bf Proof.} 
Statements (i) and (ii) follow from 5.1.
Statement (iii) from (1.4), and (iv) from
(1.6), (4.1), (4.2) and (1.1). Finally (v)
follows from the fact that $H=H_r+dr^2$. This proves the proposition.\vspace{.1in}

We are now ready to prove Theorem 1.\vspace{.2in}

\noindent {\bf Proof of Theorem 1.}\vspace{.1in}

Recall we are identifying $M-\{o\}$ with $\bS^{n-1}
\times \R^+$, and we are assuming 
$d\geq4+4\xi$.
Fix a finite atlas $\cA\0{\bS^{n-1}}$ for $\bS^{n-1}$  (see Remark 2.1(1)). 
Write $H=\cH\0{r\0{0},d}\,g\,$ which is hyperbolic on $\barB_{r\0{0}-d}$ by 5.2 (i).
Hence it remains to prove that $h$ is
$\eta$-close to hyperbolic outside
$\barB_{r\0{0}-d-1-\xi}$.\vspace{.1in}

Fix a constant $c$ such that
the metric $\hat{g}\0{r\0{0}}$ is $c$-bounded with respect to
$\cA\0{\bS^{n-1}}$. Let $p\in \R^{n}-\barB_{r\0{0}-d-1-\xi}$ and write $p=(x,r)\in
\bS^{n-1}\times \R^+$. 
Hence $r\geq r\0{0}-d-\xi-1$.
We have two cases.\vspace{.1in}

\noindent {\bf First case. $r>r\0{0}-1-\xi$.} \\
By Theorem 4.3 the metric $\cW\0{r\0{0}}g$
is radially $\eta'$-close to hyperbolic outside $\barB\0{r\0{0}-1-\xi}$  with charts 
of excess $\xi-1$,
provided
\begin{equation*}
\eta'\geq e^{27+12\xi}\,\Big(\,  e^{-2r\0{0}}\,+\, \epsilon\,\Big)
\tag{1}
\end{equation*}

\noindent In particular we can find a radially $\eta'$-close to 
hyperbolic chart $\phi$ for
$\cW\0{r\0{0}}g$ centered at $p$, with excess $\xi-1$.
But by item (iii) of 5.2 and the hypothesis $4+4\xi\leq d$
we get that $h=\cW\0{r\0{0}}\big( g\0{r\0{0}}\big)$
outside $B\0{r\0{0}-2-2\xi}$.
Therefore $\phi$ is also
an $\eta'$-hyperbolic chart for $H$ centered at $p$, with excess $\xi-1$ (see (1.8)).\vspace{.1in}

\noindent {\bf Second case. $r\leq r\0{0}-1-\xi$.} \\
Since $\hg\0{r\0{0}}$ is $c$-bounded we have that,
by Theorem 3.4, and in view of (4.1), the metric
$\cT_{_{(r\0{0}-d),d}}\bar{g}\0{r\0{0}}$
is radially $\eta''$-close to hyperbolic outside $\barB_{r\0{0}-d-1-\xi}$ with charts of excess $\xi$, provided
\begin{equation*}
\eta''\geq C\,\Big(\,  e^{-(r\0{0}-d)}\,+\, \frac{1}{d}\,\Big)
\tag{2}
\end{equation*}
\noindent Here $C=C(c,n,\xi)$.
 In particular we can find a radially $\eta''$-close to hyperbolic chart $\phi'$ for
$\cT_{_{(r\0{0}-d),d}}\bar{g}\0{r\0{0}}$ centered at $p$, with excess $\xi$.
But by Theorem 5.2(iv)
we have that $H=\cT_{_{(r\0{0}-d),d}}\bar{g}\0{r\0{0}}$
on $\barB\0{r\0{0}}$.
This together with the assumption $r\leq r\0{0}-1-\xi$ imply that $\phi'$ is also a
an $\eta''$-hyperbolic chart for $H$ centered at $p$, with excess $\xi$ (see (1.8)).
This concludes our second case.\vspace{.1in}

 Hence $H$ is radially $\eta$-close to hyperbolic with charts of excess
$\xi-1$ provided $\eta\geq \eta',\, \eta''$. Thus from (1) and (2) we 
see that we can take 
\begin{center}$\eta\geq
C\,\Big(\,e^{-(r\0{0}-d)}+\,\frac{1}{d}\,\Big)\,+\,
e^{27+12\xi}\,\Big(\,e^{-2r\0{0}}\,+\,\epsilon \Big)$\end{center}
\noindent were $C=C(c,n, \xi)$. 
But $e^{-(r\0{0}-d)}\geq e^{-2r\0{0}}$, hence the
term $e^{27+12\xi}e^{-2r\0{0}}$ in the second summand
can be absorbed into the first summand and we can
take
\begin{center}$\eta\geq
C_1\,\Big(\,e^{-(r\0{0}-d)}+\,\frac{1}{d}\,\Big)\,+\,
C_2\,\epsilon $\end{center}
\noindent with $C_2=e^{27+12\xi}$ and $C_1=C+C_2$.
This completes the proof of Theorem 1.
\vspace{.3in}

\noindent {\bf \large  Section 6. Cut Limits and Limit Metrics.}

Let $(M^{n}, g)$ be a complete Riemannian manifold with center $o\in M$. Recall that we
can write the metric on $M-\{ o\}=\R^{n}-\{0\}=\bS^{n-1}\times\R^+$ as \,\,$g=g\0{r}+dr^2$,
where $r$ is the distance to $o$.\vspace{.1in}

Fix $\xi>0$, and let $\lambda_0>1+\xi$.
Recall that the collection $\{ g_\ssl\}_{\lambda\geq   \lambda_0}$ is a 
$\odot$-family of metrics if each $g_\ssl$ is 
a  metric  of the form $g_\ssl=\big(g_\ssl\big)_r+dr^2$ defined (at least) 
for $r>\lambda-1-\xi$.\vspace{.1in}

\noindent {\bf Remark.} We will always assume that
the family
of metrics $\{ g_\ssl\}$ is smooth, that is, the map
$(x,\lambda)\mapsto g_\ssl(x)$ is smooth, $x\in M$,
$\lambda\geq\lambda_0$.
\vspace{.1in}

Recall that the $\{g_\ssl\}$ has cut limit at $b$ if there is a 
$C^2$ metric $\hg^b_{_{\infty}}$ on $\bS^{n-1}$ such that \vspace{.1in}

\noindent {\bf (6.1)}\hspace{1.5in}$
\big|\,\,{\widehat{\big(g_\ssl\big)}}_{_{\lambda+b}} \,\,-
\,\, \hg_{_{\infty}}^b\,\,\big|\,\,\longrightarrow\,\, 0
$  \,\,\,\,\,\,\,as\,\,\,\,\,\, $\lambda\ra \infty$\vspace{.1in}

\noindent {\bf Remarks.}

\noindent {\bf 1.} The metric ${\widehat{\big(g_\ssl\big)}}_{_{\lambda+b}}$ 
is the normalized spherical cut of $g_\ssl$ at
$\lambda+b$. See Section 4 or (1.4).

\noindent {\bf 2.} The arrow above means convergence in the $C^2$-norm on the 
space of $C^2$ metrics on  $\bS^{n-1}$.
See Remark 2.1(1).

\noindent {\bf 3.} The definition above implies that $\big(g_\ssl\big)_{_{r+b}}$ 
is defined for large $\lambda$,
even if $b<-1-\xi$.

\noindent {\bf 4.} Note that the concept of cut limit at $b$ depends strongly on 
the indexation of the family.\vspace{.1in}

Let $I\sbs\R$ be an interval (compact or noncompact). We say the  $\odot$-family 
$\{g_\ssl\}$ has {\it cut limits on $I$} if the convergence in (6.1) is uniform
in $b\in I$. Explicitly this means: for every $\epsilon>0$,
and $b\in I$ there are $\lambda_*$ and neighborhood
$U$ of $b$ in $I$ such that
{\scriptsize $\big|\,\,{\widehat{\big(g_\ssl\big)}}_{_{\lambda+b'}} \,\,-\,\, 
\hg_{_{\infty}}^{b'}\,\,\big|<\epsilon$}, for $\lambda>\lambda_*$ and $b'\in U$.
In particular $\{g_\ssl\}$ has a cut limit at $b$, for every $b\in I$. \vspace{.1in}

\vspace{.1in}

Consider the  $\odot$-family $\{g_\ssl\}$ and let $d>0$. Apply hyperbolic forcing to
get \begin{center} $H_\ssl=\cH\0{\lambda, d} g_\ssl $\end{center} We say that the 
family $\{H_\ssl\}$ is the {\it hyperbolically forced family} 
corresponding to the $\odot$-family $\{g_\ssl\}$.
Note that we can write $H_\ssl=( H_\ssl )_r+dr^2$.
Using Proposition 5.1 we can explicitly describe $( H_\ssl )_r$:\vspace{.1in}

\noindent {\bf (6.2.)}\,\,\,\,\,\,\,\,{\small $
\big( H_\ssl \big)_r\,\, =\,\,\left\{
\begin{array}{lll}
\big( g_\ssl \big)_r&&   \lambda +\frac{1}{2}\leq r\\
\big(1-\rho_\ssl (r)\big) \,\sinh^2r\,\, {\widehat{\big(g_\ssl \big)}}_\ssl \, +
\, \rho_\ssl(r)\,\big( g_\ssl \big)_r&
&\lambda \leq r\leq \lambda +\frac{1}{2}\\
\sinh^2r\,\, \bigg( \big( 1-  \rho_{_{(\lambda -d), d}}(r) \big)\, \sigma\0{\bS^{n-1}} \,+\,
\rho_{_{(\lambda -d), d}}(r)  \,{\widehat{\big(g_\ssl \big)}}_\ssl\,
 \bigg)&&  \lambda -d \leq r\leq \lambda  \\
\sinh^2r\,\, \sigma\0{\bS^{n-1}}&&r\leq \lambda-d
\end{array}
\right.
$}\vspace{.1in}

\noindent {\bf Proposition 6.3.}  {\it The metrics $H_\ssl$ have the following properties.

\noindent (i)  The metrics $H_\ssl$ are canonically hyperbolic on $\barB_{_{\lambda-d}}$, i.e equal 
to $\sinh^2r\,\,\sigma\0{\bS^{n-1}}+dr^2$ on
$\barB\0{\lambda-d}$, provided $\lambda >d$.

\noindent (ii) We have that $g_\ssl=H_\ssl$ outside $B_{_{\lambda+\frac{1}{2}}}.$

\noindent (iii)  The metric $H_\ssl$ coincides with $\cW\0{\lambda} \big(g\0{\lambda}\big)$ 
outside $B\0{\lambda-\frac{d}{2}}$.

\noindent (iv)  The metric $H_\ssl$ coincides with 
$\cT_{_{(\lambda-d),d}}\Big( {\overline{(g_\ssl )}}_\ssl\,\Big)$
 on $\barB\0{\lambda}$.

\noindent (v) 
If the $\odot$-family $\big\{ g_\ssl\big\}$ has  cut limit at $b=0$ then
$\big\{ H_\ssl\big\}$ has cut limits on $(-\infty,0]$. In fact we have
\begin{center}$
\hH\0{_{\infty}}^b \,\, =\,\,\left\{
\begin{array}{lll}
\hg_{_{\infty}}^0&&  b=0\\
\big( 1-  \rho(2+\frac{2b}{d})  \big)\, \sigma\0{\bS^{n-1}} 
\,+\, \rho(2+\frac{2b}{d})  \,\hg_{_{\infty}}^0 &&   -d \leq b\leq 0  \\
 \sigma\0{\bS^{n-1}}&&b\leq -d
\end{array}
\right.
$\end{center}
\noindent where $\rho$ is as in Section 3.

\noindent (vi)   If we additionally assume that $\{g_\ssl\}$ has  cut limits on 
$[0,\frac{1}{2}]$, then
$\big\{ H_\ssl\big\}$ has also  cut limits on $[0,\frac{1}{2}]$. In fact, for 
$b\in [0,\frac{1}{2}]$ we have
\begin{center}$
\hH\0{_{\infty}}^b \,\, =\,\,\big( 1-\rho(2b) \big)\,\hg_{_{\infty}}^0\,\,+\,\, \rho(2b) \,
\hg_{_{\infty}}^b
$\end{center}
\noindent where $\rho$ is as in Section 3. If $\{g_\ssl\}$ has a cut limit at 
$b> \frac{1}{2}$ then $\{H_\ssl\}$ 
has the same cut limit at $b$.

\noindent (vii)  All the rays $r\mapsto ru$, $u\in \bS^n$, emanating from the origin 
are geodesics of
$(M , H_\ssl)$. Hence, all spaces  $(M,H_\ssl)$ have center $o\in M$
and have the same geodesic rays emanating from the common center $o$. Moreover
the function $r$ (distance to  $o\in M$) is the same on all spaces $(M,H_\ssl)$.}\vspace{.1in}

\noindent {\bf Proof.}
Items (i), (ii), (iii) and (iv) follow from 
(i), (ii), (iii) and (iv) of 5.1, respectively.
Also (vii) from (v) of 5.2. We prove (v). 
We have $\hg\0{_{\infty}}^b =
lim_{\lambda\ra\infty}
{\widehat{\big(g_\ssl\big)}}_{_{\lambda+b}}$.
Let $b\leq 0$, thus $\lambda+b\leq \lambda$.
By taking $r=\lambda+b$ in
6.2 (the two cases where $r\leq \lambda
$) we get \vspace{.1in}

\noindent {\bf (6.4)}\hspace{1in}$\hH\0{_{\infty}}^b =
lim_{\lambda\ra\infty}
{\widehat{\big(H_\ssl\big)}}_{_{\lambda+b}}$

\hspace{1.41in}$=
\lim_{\lambda\ra\infty}
 \Big(  \big( 1-  \rho_{_{(\lambda -d), d}}(\lambda+b) \big)\, \sigma\0{\bS^{n-1}}\,\,+\,\,
\rho_{_{(\lambda -d), d}}(\lambda+b)  \,{\widehat{\big(g_\ssl \big)}}_\ssl  \Big)$\vspace{.1in}

\noindent But from the definition of $\rho\0{a,d}$ given
in Section 3 we have $$\rho_{_{(\lambda -d), d}}(\lambda+b)=\rho\big( 2\frac{(\lambda+b)-
(\lambda-d)}{d}\big)=\rho\big(2+\frac{2b}{d}\big)$$
\noindent Using this and taking the limit in (6.4)
we obtain (v).
\vspace{.1in}

The proof of (vi) is similar to the proof of (v).
Just note that in this case we use (6.2) in the 
case $\lambda\leq r \leq \lambda+1/2$ and
the equality $\rho\0{\lambda}(\lambda+b)=\rho(2b)$,
where $\rho$ and $\rho\0{\lambda}$ are as in Section 4.
This proves Proposition 6.3.\vspace{.2in}

\noindent {\bf Lemma 6.5.} {\it If
an  $\odot$-family $\{ g_\ssl\}$ has cut limits at $b$,
then the family
$\{\,\widehat{(g_\ssl)}\0{\lambda+b}\,\}_\ssl$
is $c$-bounded, for some $c$.}
\vspace{.1in}

\noindent {\bf Proof.} Since $\{ g_\ssl\}$ has cut limit at $b$, we have\vspace{.1in}

\noindent{\bf (6.6)}\hspace{2.4in}$\widehat{(g\0{\lambda})}\0{\lambda+b}
\ra\hat{g}\0{\infty}^b$\vspace{.1in}

Write $f_\ssl=\widehat{(g\0{\lambda})}\0{\lambda+b}$. 
Let $c\0{\infty}^b$ be such that $\hg\0{\infty}^b$ is 
$c\0{\infty}^b$-bounded. By (6.6) there is
$\lambda_*$ such that $|f_\ssl-g\0{\infty}^b|<1$,
for every $\lambda\geq\lambda_*$. Therefore
$$|f_\ssl|\leq |f_\ssl-g\0{\infty}^b|+|g\0{\infty}^b|<1+c\0{\infty}^b$$
\noindent for every $\lambda\geq\lambda_*$.
Hence  every $f_\ssl$ is $(1+c\0{\infty}^b)$-bounded, for
every $\lambda\geq\lambda_*$.
On the other hand since $\bS^{n-1}$ and the interval $[\lambda_0,\lambda_*]$ are compact, and
the map $(x,\lambda)\mapsto f_\ssl(x)$ is smooth
on $\bS^{n-1}\times [\lambda_0,\lambda_*]$,
there is $c\0{*}$ such that 
every $f_\ssl$ is $c\0{*}$-bounded, for
every $\lambda\leq\lambda_*$.
Now just take $c=max\{1+ c\0{\infty}^b, c\0{*}\}$. This proves the lemma.\vspace{.1in}

\noindent {\bf Proof of Theorem 2.}

By hypothesis the $\odot$-family $\{ g_\ssl\}$ has cut limits at $b=0$. 
This together with lemma 6.5
(take $b=0$)  imply
that there is $c$ such that
every $g_\ssl$ is $c$-bounded. We now apply
Theorem 1 to each $g_\ssl$ using the same $c$
for all $\lambda$. This proves Theorem 2.

\vspace{.1in}

Finally we give a corollary of the proof of Lemma
6.5. It is an ``interval version" of 6.5.
The proof is similar, but we have to take
track of the (now) variable $b$.\vspace{.1in}

\noindent {\bf Corollary 6.7.} {\it If
an  $\odot$-family $\{ g_\ssl\}$ has cut limits on
the compact interval $I$,
then the families
$\{\,\widehat{(g_\ssl)}\0{\lambda+b}\,\}\0{\lambda}$, $b\in I$,
are $c$-bounded for some $c$.}
\vspace{.1in}


Pedro Ontaneda

SUNY, Binghamton, N.Y., 13902, U.S.A.


{\small


\begin{thebibliography}{99}  

\bibitem{ChD} R. M. Charney and M. W. Davis, 
{\em Strict hyperbolization},
Topology {\bf 34} (1995), 329-350. 

\bibitem{O} P. Ontaneda, 
{\em Riemannian Hyperbolization}.
Arxiv: 1406.1730.

\bibitem{O1} P. Ontaneda, 
{\em On the Farrell and Jones Warping Deformation},
Journal of the London Mathematical Society {\bf 92} (2015), 566-582.


\bibitem{O2} P. Ontaneda, 
{\em Deforming an $\epsilon$-close to hyperbolic metric to a 
warped product}.
Arxiv: 1406.1741.


\bibitem{O3} P. Ontaneda, 
{\em Hyperbolic Extensions and Metrics $\epsilon$-Close to Hyperbolic}. 
Arxiv: 1406.1740. To appear in Indiana University Math. Journal.


\end{thebibliography}}
\end{document}